\documentclass[preprint,12pt]{article}
\usepackage[utf8]{inputenc}
\usepackage{indentfirst,csquotes}

\usepackage{authblk}

\usepackage{subfigure}
\usepackage{amsmath}
\usepackage{amsthm}
\usepackage{amssymb}
\newtheorem{theorem}{Theorem}[section]
\newtheorem{lemma}[theorem]{Lemma}
\newtheorem{corollary}[theorem]{Corollary}

\newtheorem{remark}[theorem]{Remark}
\newtheorem{definition}[theorem]{Definition}

\usepackage{tikz}

\usepackage{pgffor} 

\usetikzlibrary{shapes.geometric}

\usepackage{intcalc}
\usetikzlibrary{automata, positioning}

\usetikzlibrary{calc}

\usepackage{mathdots}
\usepackage{tikz-layers}

\usepackage[thinlines]{easytable}

\usepackage{float}

\usepackage{color}

\newcommand{\D}{\rm{D}}
\newcommand{\Ext}{{\rm Ext}}

\usepackage[figuresright]{rotating}

\usepackage[colorlinks]{hyperref} 
\usepackage[colorinlistoftodos]{todonotes}

\def\k{{\{k\}}}

\newtheorem{proposition}[theorem]{Proposition}
  \newtheorem{example}[theorem]{Example}

\title{Total $\{k\}$-domination in caterpillars: Structural properties and exact results} \author[1,2]{Mariana Escalante} \author[1,2]{Paola Tolomei}\author[1]{María Inés Lopez Pujato}

\affil[1]{\itshape Universidad Nacional de Rosario, Argentina}
\affil[2]{\itshape CONICET, Argentina}

 \date{}

\begin{document}

\maketitle \begingroup \renewcommand{\thefootnote}{} \footnotetext{ \textit{Email addresses:} mariana@fceia.unr.edu.ar (Mariana Escalante), ptolomei@fceia.unr.edu.ar (Paola Tolomei), lpujato@fceia.unr.edu.ar (María Inés Lopez Pujato) }
\endgroup

\begin{abstract}
The total $\{k\}$-domination problem ($T\{k\}DP$) seeks a minimum-weight function such that the sum of weights in every open neighborhood is at least $k$. While this problem is known to be NP-hard for general graphs, we demonstrate that for the family of caterpillars, it simplifies to the classical total domination number ($k=1$) via the identity $\gamma_t^{\{k\}}(G) = k\gamma_t(G)$. This reduction is supported by a general upper bound established for all graphs, which we prove to be tight specifically for caterpillars. The central tool for our analysis is the introduction of a $(T,R)$-representation—a one-to-one mapping that abstracts caterpillar structure into discrete vectors based on support vertex distribution and central path segment lengths. By characterizing these segments through congruence modulo 4, we derive an exact closed-form expression for the total $\{k\}$-domination number of any caterpillar graph.
\end{abstract}



\maketitle

\section{Introduction}\label{sec: intro}
Throughout this paper, we consider a (simple) graph $G$, with vertex set $V$ and edge set $E$. When we need to emphasize their relationship with the graph, we write $V(G)$ and $E(G)$. For each $v\in V$, $N_G(v) =\{ u\in V:\; uv\in E\}$
denotes the (open) \emph{neighborhood} of $v$ in the graph $G$, or
simply $N (v)$ when the graph is clear from the context. 
Two vertices of $G$ are \emph{adjacent} if there is an edge in $E$ having them as endpoints.
The \emph{degree} of a vertex $v$ is $d(v)=|N(v)|$. A vertex $v$ is \emph{isolated} if $d(v)=0$ and it is a \emph{leaf} if $d(v)=1$.

A (simple) \emph{path} is a graph defined by a sequence of (distinct) adjacent 
vertices. Every vertex in a path has degree 2 except for the two terminal ones, which have degree 1.  

For any $S\subseteq V$, the subgraph of $G$ \emph{induced} by $S$ is the graph with vertex set $S$ such that two vertices of $S$ are adjacent if and only if they are adjacent in $G$. 
 
A \emph{caterpillar} is a connected graph with the property that after deleting the leaves, it results in a path, called the \emph{spine} or \emph{central path} of the caterpillar. In a caterpillar $G$, every vertex in the central path adjacent to a leaf is a \emph{support vertex}.

Graph domination is a fundamental concept widely used to model a variety of problems in Combinatorial Optimization. 
Total domination in graphs was introduced by Cockayne et al. in \cite{Cockayne1980} and is now well studied in graph theory. 
A survey on this topic can be found in the so-called domination books by Haynes et al. \cite{Haynes1998, Haynes1998S}.

The study of domination and its variants is motivated by a wide range of applications in real-world network problems, including facility location, social network analysis, and resource allocation. In these contexts, a dominating set represents a collection of nodes that can monitor, serve, or influence the entire network. Total domination, in particular, is crucial for systems requiring redundancy or continuous connectivity, as it ensures that every node—including those in the dominating set—is adjacent to at least one "active" or "service" node.

While finding the minimum total dominating set is NP-hard for general graphs, the structural simplicity of caterpillar graphs enables efficient algorithmic solutions and precise combinatorial bounds. Caterpillars often serve as underlying models for linear communication backbones with peripheral nodes, such as a bus topology in local area networks or certain molecular structures in chemistry. Consequently, providing closed-form results or specialized algorithms for total domination in caterpillars not only advances theoretical graph theory but also provides tools for optimizing tree-like hierarchical architectures.

A \emph{total dominating set} of a graph $G$ (with no isolated vertices) is a set $\D$ of vertices of $G$ such that every vertex is adjacent to a
vertex in $\D$. 
The \emph{total domination problem} consists of determining the minimum cardinality of a total dominating set of $G$, denoted by $\gamma_t(G)$.

Given a positive integer number $k$, a natural variation of this concept is the total $\{k \}$-dominating function, generalizing the total dominating set when $k=1$. More precisely, given  $G$ and $k$, a function $f : V\longrightarrow \{0,1,\dots,k\}$ 
is called \emph{total  $\{k \}$-dominating function} of $G$ if for each $v\in V$, 
$f(N(v))=\sum_{u\in N(v)}f(u)\geq k.$
The value $f(u)$ is the \emph{weight} of $f$ on $u$ and for $S\subseteq V$, $f(S) = \sum_{v\in S}f(v)$. The weight of a total $\{k\}$-dominating function $f$ is $f(V)$.
The \emph{total $\{k \}$-domination number} of $G$, $\gamma_t^{\{k \}}(G)$, is the minimum weight of any total $\{k \}$-dominating function of $G$, see \cite{Li2009}.
Observe that there is a one-to-one correspondence between the vertices of $G$ at value one in a total $\{1\}$-dominating function and the vertices in a total dominating set. Then,   $\gamma_t^{\{1 \}}(G)=\gamma_t(G)$. 

Given $G$ and a  positive integer number $k$, the \emph{total $\{k\}$-domination problem} (T$\{k\}$DP) consists on finding $\gamma_t^{\{k \}}(G)$. This parameter has been studied by Aram and Sheikholeslam \cite{Aram2013}, where they provide some tight bounds for this number and determine it exactly for complete multipartite graphs. Also, Haisheng et al. in \cite{Haisheng2018} study the total $\{k\}$-domination number for other basic graph families such as cycles, paths, wheels, and pans. 
He and Liang in \cite{He2011} investigate the algorithmic complexity of computing and approximating this number. 
Regarding the computational complexity of the problem, it was shown in \cite{He2011} that determining $\gamma_t^{\{k \}}(G)$ is NP-hard for general graphs, even when restricted to certain subclasses such as chordal graphs or bipartite graphs. This complexity remains a significant hurdle for practical applications in large-scale networks, where finding optimal $\{k\}$-dominating configurations is essential for multi-node redundancy and robust coverage.

Given the problem's intractability in the general case, researchers have focused on identifying graph families for which the T$\{k\}$DP can be solved in polynomial time. Trees, and specifically caterpillar graphs, represent an important frontier in this regard. While caterpillars possess a relatively simple structure—consisting of a central spine with pendant leaves—their asymmetric nature and the flexibility of the $\{k\}$-dominating function weights $\{0, 1, \dots, k\}$ introduce non-trivial combinatorial challenges. 

The study of caterpillars in this context is not merely a theoretical exercise; it provides a foundational step for developing dynamic programming approaches and understanding how the "backbone" structure of a network influences the distribution of dominating weights. In this paper, we take advantage of the structural properties of caterpillars to determine $\gamma_t^{\{k \}}(G)$, thereby extending the existing knowledge on $\{k\}$-domination in trees.

In \cite{lagos_proc} we establish a general bound for the total $\k$-dominating number and a necessary condition for it to be tight in any graph $G$.  

\begin{proposition}[\cite{lagos_proc}]\label{prop1 y 2}
Given a graph $G$ and a positive integer number $k$, it holds that  $\gamma_t^{\{k\}}(G) \leq k\gamma_{t}(G).$ In addition, if $V$ can be partitioned into open neighborhoods of vertices of $G$ then  $\gamma_t^{\{k\}}(G)=k\gamma_{t}(G).$
\end{proposition}

\begin{figure}[h]
    \centering
\includegraphics[scale=0.25]{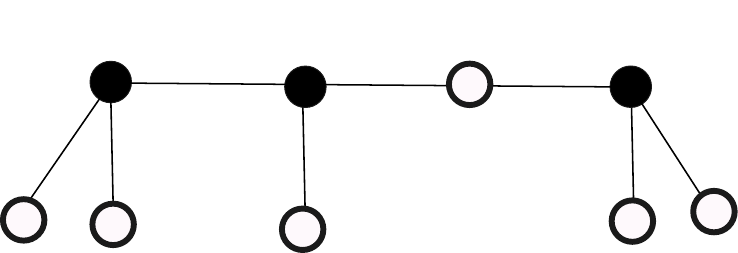}
\caption{An example of a graph where the bound in Proposition \ref{prop1 y 2} is tight.}\label{fig:disjuntas}
\end{figure}

\begin{example}\label{ejemplo}
Let us consider the graph in Figure \ref{fig:disjuntas} and a fix value $k$. We can observe that the vertices cannot be partitioned into disjoint neighborhoods. However, the bound stated in Proposition \ref{prop1 y 2} is achieved if we assign value $k$ to the black vertices and zero elsewhere. This proves that this necessary condition is not a sufficient one. 
\end{example}

It is an interesting problem to determine families of graphs that satisfy this last equality. In this line, the main results in this paper focus on the family of caterpillar graphs.

The remainder of this paper is organized as follows. In Section \ref{sec: intro}, we 
introduce the fundamental concepts of total $\{k\}$-domination and establish a 
general upper bound for $\gamma_{\{k\}t}(G)$. Section \ref{sec: D cater} focuses 
on the structural properties of $D$-caterpillars, demonstrating that for this 
family, the total $\{k\}$-domination problem can be reduced to the classical 
total domination problem ($k=1$). Section \ref{subsec: representacion} constitutes 
the central methodological framework, introducing the $(T,R)$-representation to 
encode the graph's structure and developing the necessary reduction tools. In Section 
\ref{sec: total dom general}, we synthesize these findings to present the exact 
closed-form expression for the total $\{k\}$-domination number applicable to 
any general caterpillar graph. Finally, Section \ref{sec: conclusions} concludes 
the paper and outlines directions for future research.

A preliminary version of some of the results in this contribution are presented (without proofs) in the conference paper \cite{lagos_proc}.

\section{Previous results on total $\{k\}$-domination in caterpillars}
\label{subsec: kdom total general}


Throughout this work, $[n]=\{1,\dots,n\}$ for any positive integer number $n$.

Given a caterpillar $G=(V,E)$ with $s$ support vertices, we denote $p_i$ a support vertex for every $i\in [s]$. We refer to \emph{consecutive support vertices} when they are indexed by consecutive numbers, i.e., $p_i$ and $p_{i+1}$ for some $i\in [s-1]$.

In our analysis, we consider general caterpillar graphs and caterpillars in which two consecutive support vertices do not share neighbours. When this is not the case, the two consecutive support vertices can share at most two vertices.

In the following result, we show that it suffices to consider support vertices with disjoint closed neighbourhoods, since the total $\{k\}$-domination number of a general caterpillar can be reduced to them.


\begin{lemma}\label{interseccionnodisjunta}
Let $G$ be a caterpillar with $s$ support vertices such that there exists a pair of consecutive support vertices with non-disjoint closed neighbourhoods. Let $i\in [s-1]$ such that  $N_G[p_i] \cap N_G[p_{i+1}]\neq \emptyset$.
If $v\in N_G[p_i] \cap N_G[p_{i+1}]$ with $v\neq p_i$, 
let $\Ext(G)$ be the caterpillar obtained after replacing the edge $p_iv$ by a path of the form $\{p_i,p_i^1, p_i^2,p_i^3,p_i^4,v\}$. 
Then $N_{\Ext(G)}[p_i] \cap N_{\Ext(G)}[p_{i+1}]=\emptyset$ and 
\begin{equation}\label{eq: gammat extG}
    \gamma_t^{\{k\}}(G) = \gamma_t^{\{k\}}(\Ext(G)) - 2k.
\end{equation}
\end{lemma}

\begin{figure}[htbp]
\centering
\subfigure[Case $1$.]{\includegraphics[width=60mm]{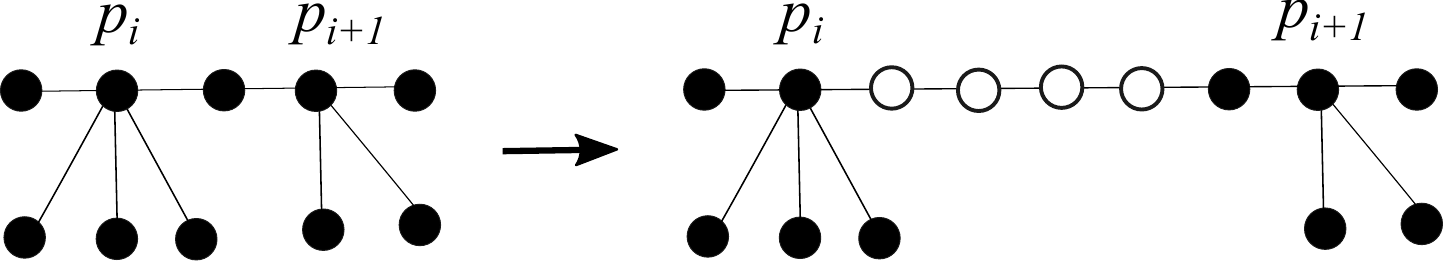}} \;\;\;\;\;\;\;\; \subfigure[Case $2$.]{\includegraphics[width=60mm]{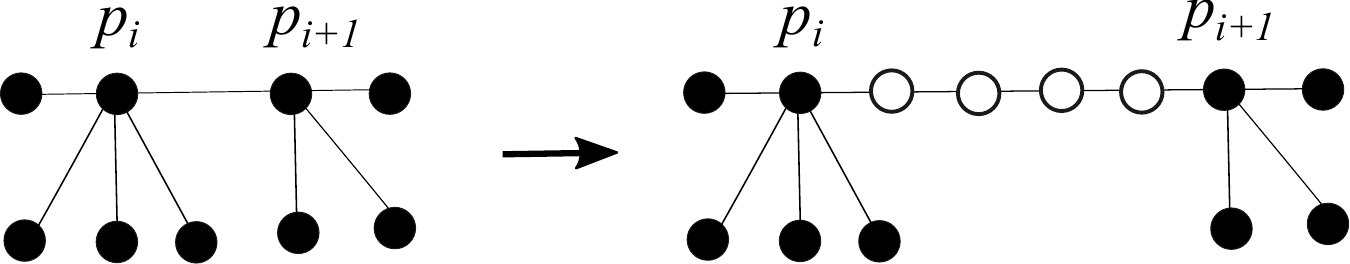}}\vspace{1mm}
\caption{Caterpillar $G$ and $\Ext(G)$ of
Lemma \ref{interseccionnodisjunta}, where the white vertices of $\Ext(G)$ correspond to $p_i^1, p_i^2,p_i^3$ and $p_i^4$, in both cases.}\label{fig 2.1}
\end{figure}

\proof
In Figure \ref{fig 2.1} we have represented the construction of $G'$ from the caterpillar $G$. We can first observe that the new caterpillar $\Ext(G)$ has two new disjoint neighbourhoods, i.e., $N(p_i^2)$ and $N(p_i^3)$. Then, 
\begin{equation}\label{eq:gama disjunto}
    \gamma_t^{\{k\}}(\Ext(G)) \geq \gamma_t^{\{k\}}(G) +2k.
\end{equation}    

Let $\D$ be an optimal total $\{k\}$-dominating set in $G$. We divide the proof into two cases, depending on the number of vertices in the intersection of the closed neighborhoods of  $p_i$ and $p_{i+1}$.

\begin{itemize}
\item Case 1: $|N_G[p_i] \cap N_G[p_{i+1}]|=1$

Let $v \in N_G[p_i] \cap N_G[p_{i+1}]$. 

If $v\in D$, consider $D'=D\cup \{p_i^1,p_i^2\}$. Clearly, $D'$ is also a total $\k$-dominating set achieving the lower bound in \ref{eq:gama disjunto}. Then, it is an optimal total $\{k\}$-dominating set of $\Ext(G)$.

If $v\notin D$, consider $D'=D\cup \{p_i^3,p_i^4\}$. Again, $D'$ is a total $\k$-dominating set achieving the lower bound in \ref{eq:gama disjunto} and then an optimal total $\{k\}$-dominating set of $\Ext(G)$. See Figure \ref{fig 2.1} $a)$.

\item Case 2: $|N[p_i] \cap N[p_{i+1}]|=2$
In this case, $N[p_i] \cap N[p_{i+1}]=\{p_{i}, p_{i+1}\}$ and $\{p_{i}, p_{i+1}\} \subseteq D$.
If we consider $D'=D\cup \{p_i^2,p_i^3\}$  then $D'$ is an optimal total $\k$-dominating set of $\Ext(G)$. See Figure \ref{fig 2.1}  $b)$.
\end{itemize}
\proof

\begin{remark}
The previous lemma can be applied iteratively to obtain a caterpillar in which every pair of consecutive support vertices has disjoint closed neighbourhoods.
\end{remark}

\begin{figure}[h]
\centering
\includegraphics[scale=0.25]{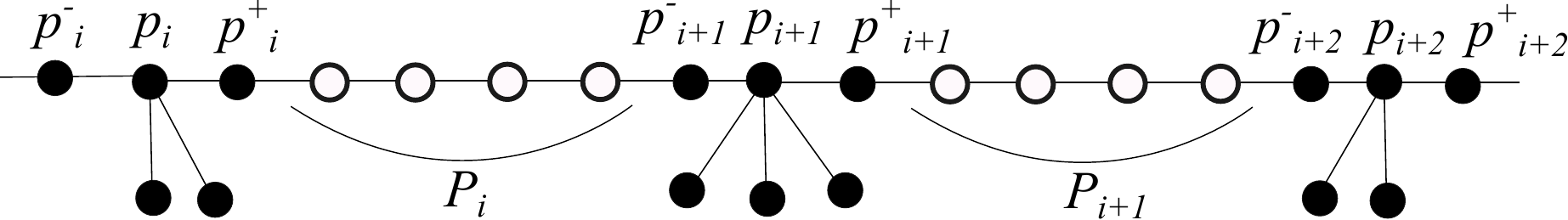}
\caption{Notation in a caterpillar.}\label{fig1}
\end{figure}

This remark allows us to work with caterpillar graphs in which the closed neighbourhoods of two consecutive support vertices are disjoint. 
Let $G$ be such a caterpillar, following the order of the vertices in the central path, for each support vertex $p_i$ with $i\in\{2,\dots,s-1\}$, let $p_i^-$ and $p_i^+$ be the neighbours of $p_i$ in $P$. Any time we use the precedence between the indices of the vertices, we consider the order of them according to the vertices in the central path starting from the vertex $p_1$ to $p_s$ (see Figure \ref{fig1}).

Under the previous notation on caterpillars, we first establish that in a caterpillar with disjoint closed support vertices neighbourhoods, there is always an optimal total $\{k\}$-dominating function having value 0 at every leaf and $k$ at every support vertex. More precisely, 

\begin{lemma}\label{lem: 1hijo}
 Given a positive integer number $k$ and a caterpillar $G=(V,E)$ such that every pair of consecutive support vertices closed neighbourhoods is disjoint, there is a total $\{k\}$-dominating function of minimum weight $f$ satisfying $f(p)=k$  and  $f(h)=0$ for every support vertex $p$ and leaf $h$ in $G$. 
\end{lemma}
\begin{proof}
Assume that there are $s$ support vertices in $G$ such that each support vertex $p_i$ has $\ell_i\geq 1$ leaves, denoted by $h_1^i,\dots,h_{\ell_i}^i$, for $i\in [s]$.
Let $g$ be a total $\{k \}$-dominating function of minimum weight in $G$. Since for every $i\in [s]$ and $j\in[\ell_i]$, $N(h^{i}_j)=\{p_i\}$,  by definition, we have that $g(p_i)=k$.

If $g(h^{i}_j)\neq 0$ for some $i\in [s]$ and $j\in [\ell_i]$ then we can define a new function $f$ on $V$ satisfying 

\[f(v)=\left\{\begin{array}{cc}
    \sum\limits_{j=1}^{\ell_i} g(h^{i}_j)+g(p_i^{+}) & \text{ if } \;v=p_i^{+} \text{for} \;i\in [s-1],\\
    \sum\limits_{j=1}^{l_s} g(h^{s}_j)+g(p_s^{-}) & \text{ if } \;v=p_s^{-}, \\
 0&  \text{ if } v=h_j^i\; \text{ for } j \in [\ell_i]\; \text{ and } \;i\in [s],\\
    g(v)& \text{ otherwise}.
\end{array}\right.
\]

    Clearly, $f$ is a total $\{k \}$-dominating function and it is optimal since $f(V)=g(V)$.
\end{proof}

\begin{remark}\label{primerremark}
  By Lemma \ref{lem: 1hijo}, the weight of a total $\{k \}$-dominating function in a caterpillar is independent of the number of leaves adjacent to a support vertex. Therefore, we can assume,  w.l.o.g., that every support vertex $p_i$ has only one leaf adjacent to it, and denote it as $h_i$ for each $i\in [s]$.  
\end{remark}

\begin{figure}[h]
\centering
\includegraphics[scale=0.24]{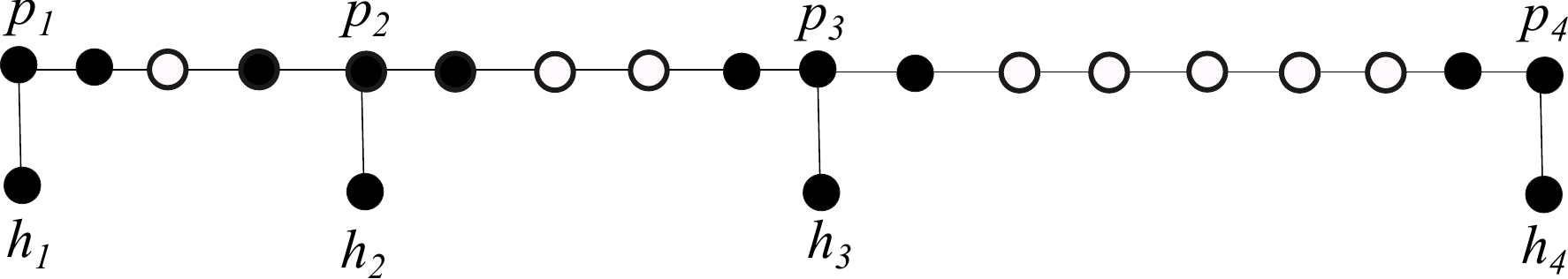}
\caption{$\D$-caterpillar.}\label{fig 2.3}
\end{figure}

\begin{definition} 
Let $G$ be a caterpillar such that every support vertex has one leaf adjacent to it. Moreover, if, for any two consecutive support vertices of $G$, their closed neighborhoods are disjoint, then $G$ is called a \emph{$\D$-caterpillar}. See Figure \ref{fig 2.3}.
\end{definition}

From now on, given a $\D$-caterpillar, we denote by $p_i$ its $i$-th support vertex and by $h_i$ the unique leaf adjacent to it.

\subsection{On total $\{k\}$-domination in $\D$-caterpillars}\label{sec: D cater}

Let us introduce some useful notation for every $\D$-caterpillar $G$.

\begin{definition}\label{def: GiPi}
Let $G=(V,E)$ be a $\D$-caterpillar with $s$ support vertices. The support vertices in the central path $P$ are denoted by $p_i$ for $i\in\{2,\dots,s-1\}$, $p_i^-$ and $p_i^+$ the neighbours of $p_i$ in $P$. Let $P_i$ be the set of vertices in $P$ strictly between $p_i^+$ and $p_{i+1}^-$, for every $i\in [s-1]$. Clearly,  when $P_i\neq \emptyset$, it induces a path in $G$ for every $i\in [s-1]$. 

Let $G_i$ be the subgraph of $G$ induced by 
\begin{equation}\label{eq: Vi}
    V_i=\{h_i,p_i,p_i^+\}\cup P_i\cup \{p_{i+1}^-, h_{i+1},p_{i+1}\} 
\end{equation} and $n_i=|V_i|$ for each $i\in [s-1]$. 
Note that each graph $G_i$ is an induced path. 
\end{definition}
\begin{figure}[h]
\centering
\includegraphics[scale=0.25]{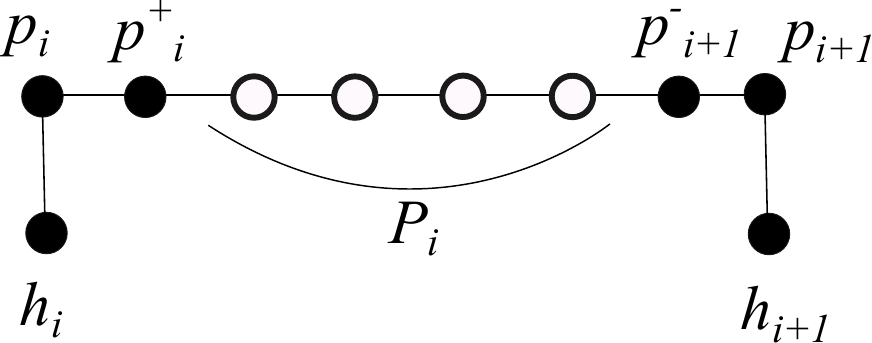}
\caption{The subgraph $G_i$ of $G$.}
\label{fig: Gi}
\end{figure}

\begin{figure}[htbp]
\centering
\includegraphics[scale=0.25]{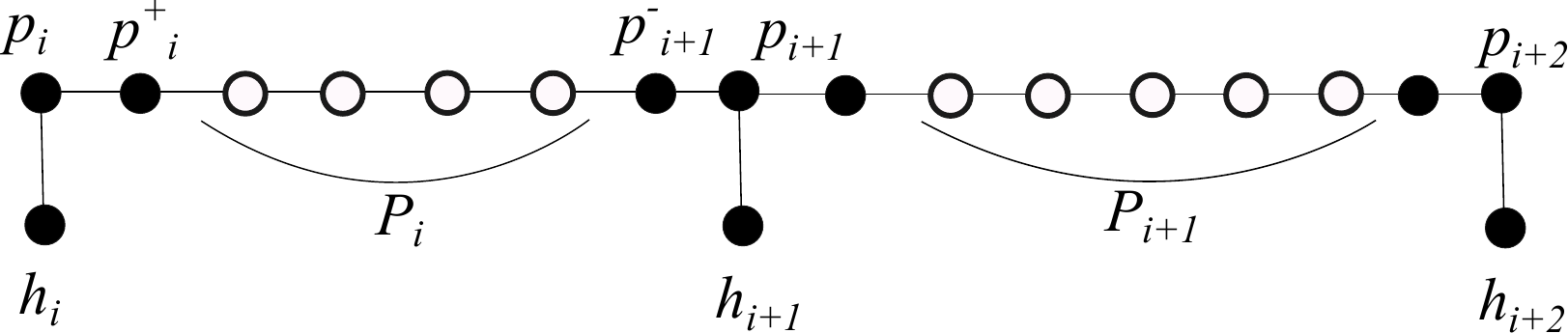}
\caption{The white nodes correspond to $P_i$ and $P_{i+1}$  in $G_i\oplus G_{i+1}$.}
\label{fig: GiGi+1}
\end{figure}

According to the previous definition, $G$ can be decomposed into $s-1$ subgraphs of the form $G_i$, where the vertices of $G_i$ and $G_{i+1}$ clearly intersect in $\{p_{i+1},h_{i+1}\}$ for every $i\in [s-2]$. To emphasize this decomposition, we write $G=\oplus_{i=1}^{s-1}G_i$.
Since $G$ is a $\D$-caterpillar, $n_i\geq 6$ for each $i\in [s-1]$. We represent a graph $G_i$ in Figure \ref{fig: Gi} and $G_i\oplus G_{i+1}$, the paths induced by $P_i$ and $P_{i+1}$ and the leaves $h_i$ and $h_{i+1}$ in Figure \ref{fig: GiGi+1}. 

When the $\D$-caterpillar has only two support vertices, the graph itself is a path with at least 6 vertices. Based on the results in \cite{Haisheng2018}, Lemma 2.5 in \cite{lagos_proc} establishes the existence of an optimal total $\{k\}$-dominating function $\{0,k\}$-valued and also a closed formula for the total $\{k\}$-dominating number of paths. 
We rewrite this result here for completeness.

 \begin{lemma} \label{lem: camino} If $G$ is a path of $n\geq 6$ vertices, then:
    \begin{enumerate}
        \item[a)]     
    there exists an optimal total $\{k\}$-dominating function that is $\{0,k\}$-valued, and

 \item[b)] $\gamma_t^{{\{k\}}}(G)=\left\{
\begin{array}{lll}

k\frac{n+2}{2}&& \text{if} \;\;n-6\equiv 0\;\; (\operatorname{mod} 4),\\
&&\\
k\frac{n+1}{2}&& \text{if} \;\;n-6\equiv 1 \;\; \text{or}\;\; n-6\equiv 3\;\;(\operatorname{mod} 4),\\
&&\\
k\frac{n}{2}&&\text{if} \;\;n-6\equiv 2\;\;(\operatorname{mod} 4).\\
\end{array}\right .$
\end{enumerate}
\end{lemma}

We extend this result to the case of more than two support vertices in the given $\D$-caterpillar. 

The following theorem was originally presented in \cite{lagos_proc}. However, we include a revised proof of it here, since it provides some intuition on how to define an optimal total $\{k\}$-dominating function, central to our subsequent results.

\begin{theorem}\label{thm: 0-k}
Let $G$ be a $\D$-caterpillar. Then, for every positive integer number $k$, there exists a $\{0,k\}$-valued optimal total $\{k\}$-dominating function.
\end{theorem}

\begin{proof}
We proceed by induction on the number of support vertices $s$ in the $\D$-caterpillar $G$. For the base case $s = 2$, the result follows directly from the fact that the caterpillar is a path. 

Assume that the statement holds for any $\D$-caterpillar with at most $j$ support vertices. Let $G$ be a $\D$-caterpillar with $s = j + 1$ support vertices, which can be viewed as the composition $G = G' \oplus G_j$, where $G'$ is a $\D$-caterpillar containing the first $j$ support vertices and $G_j$ is the subgraph induced by $V_j$ associated with the last two support vertices $p_j$ and $p_{j+1}$. 

By the inductive hypothesis, there exists an optimal total $\{k\}$-dominating function $g$ for $G'$ such that $g(v) \in \{0, k\}$ for all $v \in V(G')$. Our goal is to extend or modify $g$ to a function $f$ on $V(G)$ that remains $\{0, k\}$-valued and optimal.

Let $a = g(p^-_j)$ be the weight assigned to the neighbor of $p_j$ in the central path of $G'$. Based on the weight $g$ assigns to the neighborhood of $p_j$, we define $f$ on $V(G') \setminus \{h_j\}$ to coincide with $g$, while setting $f(h_j) = 0=f(h_{j+1})$ to maintain the $\{0, k\}$-property on $G$. To determine the values of $f$ on the path $P_j$ (the vertices between $p^+_j$ and $p^-_{j+1}$) and the final support vertex $p_{j+1}$, let $|P_j| = 4t_j + r_j$ with $r_j \in \{0, 1, 2, 3\}$. 

The total weight of $f$ on $V_j \cup \{p^-_j\}$ must satisfy the following lower bound arising from the total $\{k\}$-domination constraints on $p_j$, $p_{j+1}$, and the $2t_j$ disjoint open neighborhoods within $P_j$:
\begin{equation}
f(V_j \cup \{p^-_j\}) \geq 4k + 2kt_j + \sum_{l=1}^{r_j} f(w_{4t_j+l})
\end{equation}
where $w_l$ are the vertices in the path induced by $P_j$.

We define $f$ on the path induced by $P_j = \{w_1, \dots, w_{4t_j+r_j}\}$ by alternating weights according to the parity and the boundary value $a \in \{0, k\}$:
$f(w_{4t+1}) = 0$, $f(w_{4t+2}) = a$, $f(w_{4t+3}) = k$, and $f(w_{4t+4}) = k-a$, for $t \in \{0, \dots, t_j-1\}$.

The definition of $f$ for the remaining vertices depends on $r_j$:
\begin{itemize}
    \item \textbf{Case $r_j = 0$:} Setting $f(p^-_{j+1}) = f(h_{j+1}) = 0$ achieves the lower bound, ensuring optimality.
    \item \textbf{Case $r_j = 1$:} If $a=k$, we set $f(p^-_{j+1})=k$ and $f(h_{j+1})=0$. If $a=0$, we set $f(p^-_{j+1})=0$ and $f(h_{j+1})=k$.
    \item \textbf{Case $r_j = 2$:} Setting $f(p^-_{j+1})=k$ and $f(h_{j+1})=0$ satisfies the constraints and reaches the lower bound.
    \item \textbf{Case $r_j = 3$:} If $a=0$, we set $f(w_{4t_j+3}) = f(p^-_{j+1}) = k$. If $a=k$, we set $f(w_{4t_j+2}) = f(w_{4t_j+3}) = k$.
\end{itemize}

In all scenarios, the resulting function $f$ is a total $\{k\}$-dominating function of $G$ with $f(v) \in \{0, k\}$ for all $v \in V(G)$. Since the construction matches the general lower bound for total $\{k\}$-domination in these structures, $f$ is optimal. This completes the induction.
\end{proof}

Note that the previous result allows us to study total domination instead of total $\{k\}$-domination in $\D$-caterpillars. More precisely, 

\begin{corollary}\label{cor: reduccion}
Given a positive integer number $k$ and a $\D$-caterpillar $G$, the T$\{k\}$DP can be reduced to the total domination problem of $G$. Moreover,  $\gamma_t^{\{k\}}(G)=k\gamma_{t}(G)$.
\end{corollary}

\begin{proof}

Let $G=(V,E)$ be a $\D$-caterpillar and let $S$ be an optimal total dominating set of $G$ of size  $L$. Define the function $f$ on $V$ as $f(v)=k$ for all $v \in S$ and $f(v)=0$ for all $v \in V \setminus S$. Then, it is clearly an optimal total $k$-dominating function of $G$ with weight $kL$.

Conversely, by Theorem $\ref{thm: 0-k}$ let $f$ be an optimal total $\{k\}$-dominating function of $G$ $\{0,k\}$-valued  weight $kL$. Then the set $S=\{v \in V : f(v)=k\}$ is clearly an optimal total $\{k\}$-dominating set of $G$ of size  $L$.
\end{proof}

The previous theorem implies that the bound given in Proposition \ref{prop1 y 2} for the total $\{k\}$-domination number is achieved for the family of all $\D$-caterpillars. Also, this result shows that, to study $\k$-domination in this family, we can restrict our attention to the problem of total domination in graphs.  

\section{Total domination in $\D$-caterpillars}\label{sec: dom total}

In the proof of Theorem \ref{thm: 0-k}, we observed that the congruence modulo 4 of the number of vertices in the paths induced between consecutive support vertices plays a decisive role in determining the optimal weight distribution. This section formalizes this observation by introducing a compact representation for $\D$-caterpillars.

\subsection{On the representation of a $\D$-caterpillar}
\label{subsec:representation}

While the previous results establish the reduction of the total $\{k\}$-domination problem to the case $k=1$, a significant challenge remains: providing an exact closed-form expression for $\gamma_t(G)$ for any caterpillar $G$. The difficulty lies in the asymmetric distribution of weights across the central path and the sensitivity of the total weight to the specific lengths of the paths between support vertices.

To address this, we introduce the $(T,R)$-representation, a combinatorial framework that abstracts the nature of a $\D$-caterpillar into two discrete vectors. 
By considering congruence modulo 4 in the segments of the central path between the open neighbourhood of two consecutive support vertices, we define the integer-valued vectors $T$ and $R$, 
and provide a one-to-one mapping between the family of $\D$-caterpillars and the specific discrete space. This representation is not merely descriptive; it serves as a fundamental tool to derive the exact total domination number.

\begin{lemma}
\label{lemma:representation}
There is a one-to-one correspondence between the family of $\D$-caterpillars with $s$ support vertices and the set $\mathbb{N}_0^{s-1} \times \{0, 1, 2, 3\}^{s-1}$.
\end{lemma}

\begin{proof}
Let $(T,R) \in \mathbb{N}_0^{s-1} \times \{0, 1, 2, 3\}^{s-1}$, where $T = (t_1, t_2, \dots, t_{s-1})$ and $R = (r_1, r_2, \dots, r_{s-1})$. We define $G(T,R)$ as the $\D$-caterpillar having $s$ support vertices $p_i$ with $i\in[s]$ where each vertex subset $P_i$ between the open neighbourhoods of $p_i$ and $p_{i+1}$ in the central path of $G(T,R)$, satisfies $|P_i| = 4t_i + r_i$ for each $i \in [s-1]$. 
This completes the definition of the $\D$-caterpillar $G(T,R)$ obtained after the two vectors $T$ and $R$. For instance, the $\D$-caterpillar depicted in Figure \ref{fig: GiGi+1} corresponds to $T=(1,1)$ and $R=(0,1)$, considering $i=1$.

Conversely, let $G$ be a $\D$-caterpillar with a central path $P$ and support vertices $p_i$, for $i\in [s]$. 
Let $P_i$ be the set of vertices in $P$  between the open neighbourhoods of $p_i$ and $p_{i+1}$, for each $i\in [s-1]$.  If $t_i\in \mathbb{N}_0$ and $r_i\in \{0,1,2,3\}$ are such that  $|P_i|= 4t_i+r_i$ for each $i\in [s-1]$, we define $T=(t_1, t_2, \dots, t_{s-1})$ and $R=(r_1, r_2, \dots, r_{s-1})$. 

This shows the one-to-one correspondence between the pair $(T,R)$ and the $\D$-caterpillar $G$.
\end{proof}

\begin{remark}
Given $G$, we refer to this unique pair $(T(G),R(G))$ as the $(T,R)$-representation of $G$. According to the previous result, given $G$ a $\D$-caterpillar we have that
$T=(t_i)_{s-1}\in \mathbb{N}_0^{s-1}$ and $R=(r_i)_{s-1}\in  \{0,1,2,3\}^{s-1}$, where $|P_i|=4t_i+r_i$ for every $i\in [s-1]$.
\end{remark}

\begin{figure}[h]
\centering
\includegraphics[scale=0.25]{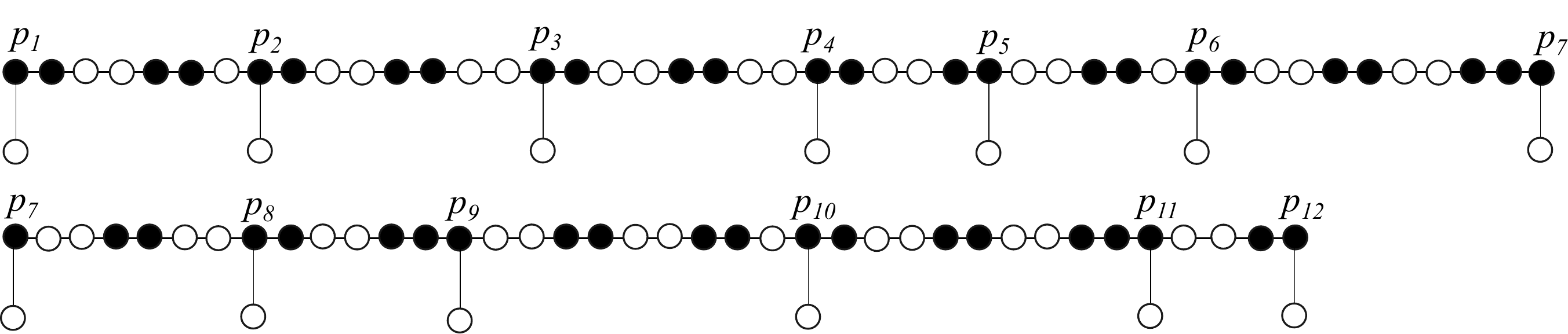}
\caption{The $\D$-caterpillar in 
Example \ref{ex:12parents}, where the black vertices correspond to an optimal total dominating set.}
\label{fig: ejemcater}
\end{figure}

\begin{definition}
\label{def:sets_rj}
Given the $(T,R)$-representation of a $\D$-caterpillar $G$, we define the sets \[\mathcal{R}_j = \{G_i : r_i = j, \text{ for } i \in [s-1]\}\]
and their corresponding index sets $I_j = \{i : G_i \in \mathcal{R}_j\}$, for each $j \in \{0, 1, 2, 3\}$. 
Furthermore, if $\mathcal{R}_3 \neq \emptyset$, we denote by $\mathcal{G}_j$ the $j$th-maximal set of $\ell_j$ consecutive graphs $G_i \in \mathcal{R}_3$, for $j \in [m]$, where $m$ indicates the number of maximal subsequences of consecutive components equal to 3 in the vector $R$.
\end{definition}
We illustrate this definition in an example; see Figure \ref{fig: ejemcater}. 
\begin{example}
\label{ex:12parents}
Consider a $\D$-caterpillar $G$ with 12 support vertices defined by the following $(T,R)$-representation:
\begin{itemize}
    \item $T = (1, 1, 1, 0, 0, 1, 1, 0, 1, 1, 0)$
    \item $R = (0, 1, 1, 2, 3, 3, 0, 3, 3, 3, 1)$
\end{itemize}
Based on the components of the remainder vector $R$, we can partition the subgraphs $G_i$ into the index sets:
\begin{itemize}
    \item $I_0 = \{1, 7\}$, $I_1 = \{2, 3, 11\}$, $I_2 = \{4\}$, and $I_3 = \{5, 6, 8, 9, 10\}$.
\end{itemize}
Furthermore, the vector $R$ contains two maximal sequences of consecutive 3s ($m=2$):
\begin{itemize}
    \item $\mathcal{G}_1 = \{G_5, G_6\}$ with length $\ell_1 = 2$.
    \item $\mathcal{G}_2 = \{G_8, G_9, G_{10}\}$ with length $\ell_2 = 3$.
\end{itemize}
\end{example}

Given a $\D$-caterpillar $G$ we define a smaller  $\D$-caterpillar $G'$ having a $\{0,2,3\}$-valued vector $R(G')$ in its $(T(G'),R(G'))$-representation. This reduction will be useful for determining the total domination number of $G$. See Theorem \ref{thm: main}.

\begin{definition}
\label{def:reduced_graph}
Given the $(T,R)$-representation of a $\D$-caterpillar $G$, we define the vector $R'$ deleting from $R$ all components at value $1$.
The vector $T'$ is defined accordingly by preserving the components of $T$ that correspond to the indices kept in $R'$. For simplicity, we denote by $G'$ the $\D$-caterpillar $G(T', R')$ and we refer to it as the reduced graph obtained from $G$.
\end{definition}

\begin{example}\label{ex:G'}
Given the $(T,R)$-representation of the $\D$-caterpillar $G$ in Example \ref{ex:12parents}, according to the previous definition, $R'=(0,2,3,3,0,3,3,3)$ and $T'=(1,0,0,1,1,0,1,1)$. Then, $G'=G_1\oplus G_4\oplus G_5\oplus G_6\oplus G_7 \oplus G_8 \oplus G_9 \oplus G_{10}$, see Figure \ref{fig: exG'}. 

\begin{figure}[h]
\centering
\includegraphics[scale=0.22]{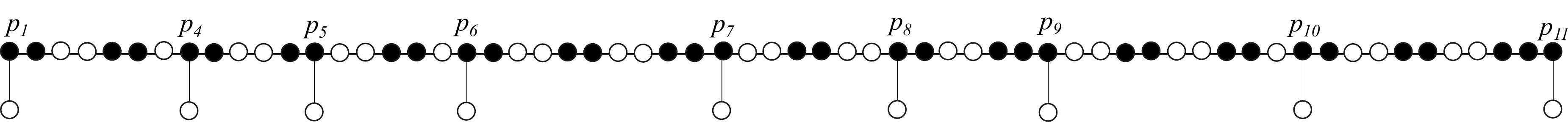}
\caption{The $\D$-caterpillar $G'$ with $(T',R')$-representation, where the black vertices correspond to an optimal total dominating set.}
\label{fig: exG'}
\end{figure}
\end{example}

\subsection{On the total domination number in $\D$-caterpillar}

We begin this section by presenting some properties satisfied by an optimal total dominating function on a $\D$-caterpillar that will be useful for presenting the main result of this section, i.e., the total domination number of a general $\D$-caterpillar. Recall that in a $\D$-caterpillar with support vertices $p_i$ for $i\in [s]$,  $p_i^-$ is \emph{previous} neighbour of $p_i$ in the central path of $G$ and $p_i^+$ is \emph{subsequent} neighbour of $p_i$ in the central path of $G$. 



According to the particular neighbourhood of support vertices in $\D$-caterpillars, in the following lemma, we show that we can consider two different possible assignments in an optimal total dominating function.

\begin{lemma}\label{primerclaim}
If $G$ is a $\D$-caterpillar with support vertices $p_1,\dots,p_s$, then there is an optimal total dominating function $f$ satisfying  either 
$f(p_i)=f(p_i^+)=1$ or $f(p_i)=f(p_i^-)=1$, for every $i\in [s-1]$.
\end{lemma}
\begin{proof}
By Lemma \ref{lem: 1hijo}, there is an optimal total dominating function $f$ such that $f(h_i)=0$ and $f(p_i)=1$ for all $i\in [s]$. Since $N(p_1)=\{p_1^+\}$, $N(p_i)=\{p_i^{-}, p_i^{+}\}$ for all $i\in \{2,\dots,s-1\}$ and $N(p_s)=\{p_s^{-}\}$ the optimality of $f$ proves the result.
\end{proof}

If the $\D$-caterpillar is decomposed according to Definition \ref{def: GiPi}, the previous lemma helps in determining the value of an optimal total dominating function on each graph $G_i$.

\begin{lemma}\label{segundoclaim}
Let $G$ be a $\D$-caterpillar with support vertices $p_1,\dots, p_s$ such that $G=\oplus_{i=1}^{s-1} G_i$ and $P_i=V(G_i)\setminus(N[p_i]\cup N[p_{i+1}])$ for $i\in [s-1]$. If  $t_i$ and $r_i$ are the $i$-th components of $T$ and $R$ in the  $(T,R)$-representation of $G$, there is an optimal total dominating function $f$ satisfying.
\begin{enumerate}
\item[i)]  
\begin{equation}\label{pesocamino1} 
f(P_i)=\left\{
\begin{array}{lll}

2t_i+1&& \text{ if }  r_i=3, \\

2t_i&& \text{otherwise};
\end{array}\right .
\end{equation}
if $f(p_i)=f(p_i^{+})=1$, and 
\item[ii)]   
\begin{equation}\label{pesocamino2} 
f(P_i)=\left\{
\begin{array}{lll}
2t_i&& \text{ if }  r_i\in \{0,1\},\\

2t_i+1&& \text{ if }  r_i=2, \\ 

2t_i+2&& \text{ if }  r_i=3,\\

\end{array}\right .
\end{equation}
if $f(p^{-}_i)=f(p_i)=1$.
\end{enumerate}

\end{lemma} 
\begin{proof}

Given the $\D$-caterpillar $G$, in its $(T,R)$-representation, the components of $T$ and $R$, $t_i$ and $r_i$, are the non-negative integer  numbers such that $|P_i|=4t_i+r_i$ and $r_i\leq 3$, for $i\in [s-1]$. Let $i\in [s-1]$ and assume that $P_i=\{w_1,\dots,w_{{b}_i}\}$.

According to Lemma \ref{primerclaim} there exists an optimal total dominating function $\tilde{f}$  satisfying either that $\tilde{f}(p_i)=\tilde{f}
(p^+_i)=1$ or $\tilde{f}(p^{-}_i)=\tilde{f}(p_i)=1$. In addition, it holds that
\begin{equation}\label{cota_camino}
\tilde{f}(P_i)\geq 2t_i+\sum_{l=1}^{r_i} \tilde{f}(w_{4t_i+l}).    
\end{equation}

\begin{enumerate}
\item[$i)$] Suppose first that  $\tilde{f}(p_i)=\tilde{f}
(p^+_i)=1$. 

If $r_i\in\{0,1\}$ we define a total dominating function $f$ as
\[
f(v)=\left\{
\begin{array}{ll}
 1 & \text{ if } v\in\{w_{4t+3},w_{4t+4}\}, \\
    0 & \text{ if }  v\in \{w_{4t+1},w_{4t+2},p_{i+1}^-\}\cup \{w_{4t_i+l}\}:l\in [r_i],\\
      \tilde{f}(v) & \text{ otherwise}.  
\end{array}
\right.
\]
for $t\in\{0,\dots, t_i-1\}$, achieves the value $f(P_i)=2t_i$, thus proving its optimality due to (\ref{cota_camino}). 

Note that $[r_i]=\emptyset$ if $r_i=0$.

In case $r_i=2$, 
\[]
f(v)=\left\{
\begin{array}{ll}
 1 & \text{ if } v\in\{w_{4t+3},w_{4t+4},p_{i+1}^-\}, \\
    0 & \text{ if }  v\in \{w_{4t+1},w_{4t+2},w_{4t_i+1},w_{4t_i+2}\},\\
  \tilde{f}(v) & \text{ otherwise}.  
\end{array}
\right.
\]
for $t\in\{0,\dots, t_i-1\}$,
achieves the value $f(P_i)=2t_i$. Again, this proves its optimality due to (\ref{cota_camino}).

Similarly, when $r_i=3$, we define the total dominating function
$f$ as
\[
f(v)=\left\{
\begin{array}{ll}
 1 & \text{ if } v\in\{w_{4t+3},w_{4t+4},w_{4t_i+3},p_{i+1}^-\}, \\
    0 & \text{ if }  v\in \{w_{4t+1},w_{4t+2},w_{4t_i+1},w_{4t_i+2}\},\\
  \tilde{f}(v) & \text{ otherwise}.  
\end{array}
\right.
\]
for $t\in\{0,\dots, t_i-1\}$,
achieves the lower bound $2t_i+1$ given by (\ref{cota_camino}), and then proving its optimality.

\item[$ii)$] Assume now that $\tilde{f}(p^{-}_i)=\tilde{f}(p_i)=1$.

If $r_i=0$, we define the total dominating function $f$ as 
\[
f(v)=\left\{
\begin{array}{ll}
 1 & \text{ if } v\in\{w_{4t+2},w_{4t+3}\}, \\
    0 & \text{ if }  v\in \{w_{4t+1},w_{4t+4},p_{i+1}^-\},\\
  \tilde{f}(v) & \text{ otherwise}.  
\end{array}
\right.
\]
for $t\in\{0,\dots, t_i-1\}$,
achieves the value $f(P_i)=2t_i$, thus proving its optimality due to (\ref{cota_camino}).
\end{enumerate}

Similarly, when $r_i=1$,  
the total dominating function $f$ is defined as
\[
f(v)=\left\{
\begin{array}{ll}
 1 & \text{ if } v\in\{w_{4t+2},w_{4t+3},p_{i+1}^-\}, \\
    0 & \text{ if }  v\in \{w_{4t+1},w_{4t+4}\},\\
  \tilde{f}(v) & \text{ otherwise}.  
\end{array}
\right.
\]
for $t\in\{0,\dots, t_i-1\}$,
and it achieves the value $f(P_i)=2t_i$ thus proving its optimality due to (\ref{cota_camino}).

When $r_i=2$, the function $f$ is defined as
\[
f(v)=\left\{
\begin{array}{ll}
 1 & \text{ if } v\in\{w_{4t+2},w_{4t+3},p_{i+1}^-,w_{4t_i+2}\}, \\
    0 & \text{ if }  v\in \{w_{4t+1},w_{4t+4},w_{4t_i+1}\},\\
  \tilde{f}(v) & \text{ otherwise}.  
\end{array}
\right.
\]
for $t\in\{0,\dots, t_i-1\}$.

Finally, when $r_i=3$ the function $f$ is defined as
\[
f(v)=\left\{
\begin{array}{ll}
 1 & \text{ if } v\in\{w_{4t+2},w_{4t+3}, w_{4t_i+2},w_{4t_i+3}\}, \\
    0 & \text{ if }  v\in \{w_{4t+1},w_{4t+4},w_{4t_i+1},p_{i+1}^-\},\\
  \tilde{f}(v) & \text{ otherwise}.  
\end{array}
\right.
\]
for $t\in\{0,\dots, t_i-1\}$.

In both cases, $f$ is clearly a total dominating function of $G$. Also note that when $r_i\in\{2,3\}$, the set $P_i\cup\{p_{i+1}^-\}$ can be partitioned into $2t_i+ (r_i-1)$ disjoint open neighbourhoods of vertices of $G$. Then, this function achieves optimality using the bound provided by (\ref{cota_camino}).
\end{proof}

In the next observation, we explicitly show that there is always an optimal total dominating function that takes into account the two possible situations considered in Lemma \ref{primerclaim} to define its value in the neighbourhood of the following support vertex. More precisely, after the proof of Lemma \ref{segundoclaim},  we can state the following.

\begin{remark}\label{rem:inicios}
There is an optimal total dominating function $f$ satisfying, for $i\in[s-1]$, 

\begin{enumerate}
\item[$i)$]
 \[
\left\{
\begin{array}{lll}
f(p_{i+1})=f(p_{i+1}^{+})=1 &&
\text{ if }  r_i=0 \text{ or } r_i=1, \\

f(p_{i+1})=f(p_{i+1}^{-})=1&& \text{otherwise};
\end{array}\right. 
\]
if\;\; $f(p_i)=f(p_i^{+})=1$.

\item[$ii)$]  
\[
\left\{
\begin{array}{lll}
f(p_{i+1})=f(p_{i+1}^{+})=1 &&
\text{ if }  r_i=0 \text{ or } r_i=3, \\
f(p_{i+1})=f(p_{i+1}^{-})=1&& \text{otherwise}.
\end{array}\right.\]
 if $f(p^{-}_i)=f(p_i)=1$.
 
\end{enumerate}
In both cases, we have considered $p_{s}^{+}=h_s$.

\end{remark}

According to the previous result we are in a position of determining the total domination number of the caterpillar for which its $(T,R)$-representation coincides with the $(T',R')$-representation given in Definition \ref{def:reduced_graph}. This provides a basic reduction for obtaining the total domination number for general caterpillar graphs.

Recall that, given a $\D$-caterpillar such that $G=\oplus_{i=1}^{s-1} G_i$, the value of the total domination number of any $G_i$ can be obtained by the results in \cite{Haisheng2018}, for each $i\in[s-1]$ since they are actually paths in the original graph $G$.
Recall that we have already know the total domination number of paths and then then we consider $s\geq 3$ in the next proposition and the particular situation where $I_0=I_1=\emptyset$. We fundament this assumption in the subsequent results.

\begin{proposition}\label{prop: G'}
Let $G$ be a $\D$-caterpillar with $s\geq 3$ support vertices, with $T=T(G)$ and $R=R(G)$ its $(T,R)$-representation and $I_0=I_1=\emptyset$.

If $I_2\neq \emptyset$,
\begin{equation}\label{eq: gamma reducido}
\gamma_t(G)=2t_1+\sum_{i=2}^{s-1} (2t_i+E_{i})+ \sum_{j=1}^{m} \left(2\ell_j - \lfloor \frac{\ell_j}{2}\rfloor\right)+2s,
\end{equation}
 where  
\[
E_{i}=\left\{
\begin{array}{ll}
 1 & \text{ if } i-1, i \in I_2,  \\
 0 & \text{ otherwise}. 
\end{array}
\right.
\]
for $i\in [2,s-1]$.

If $I_2=\emptyset$, 
\begin{equation}\label{eq: gamma reducido sin I2}
\gamma_t(G)=\sum_{i=1}^{s-1}2 t_i+\lfloor\frac{s-1}{2}\rfloor-1 + 3s.
\end{equation}

\end{proposition}

\begin{proof} 
Lemma \ref{primerclaim} and Lemma \ref{segundoclaim} imply that there exists an optimal total dominating function $f$ such that: $f(p_1)=f(p_1^{+})=1$ and 

\begin{equation}\label{eq: gamma prop}
    \gamma_t(G)=\sum_{i=1}^{s-1} f(P_i)+2s.
\end{equation}\label{eq:gamma}

Assume first that $I_2\neq \emptyset$.

By hypothesis, we have to consider only four combinations of $r_{i-1}$ and $r_{i}$, for $i\in[2,s-1]$, i.e., $r_{i-1}\in\{2,3\}$ and $r_{i}\in\{2,3\}$. 

\begin{itemize}

\item[-] Case $r_{i-1}=2$ and $r_{i}=2$.

From Remark \ref{rem:inicios}, $f(p^{-}_{i})=f(p_{i})=1$ and then,  Lemma\ref{segundoclaim} implies that $f(P_{i}) = 2t_{i}+1=2t_{i}+E_{i}$.

\item[-] Case $r_{i-1}=2$ and $r_{i}=3$. 

From Remark \ref{rem:inicios}, $f(p^{-}_{i})=f(p_{i})=1$ and then  Lemma\ref{segundoclaim} implies that $f(P_{i}) = 2t_{i}+2=2t_{i}+E_i+2$. 


\item[-] Case $r_{i-1}=3$ and $r_{i}=2$. 

If $f(p^{-}_{i-1})=f(p_{i-1})=1$ then $f(p^{+}_{i})=f(p_{i})=1$, and  Lemma \ref{segundoclaim} implies that $f(P_{i}) = 2t_{i}$.

If $f(p^{+}_{i-1})=f(p_{i-1})=1$ then  $f(p^{-}_{i})=f(p_{i})=1$, implying that $f(P_{i}) = 2t_{i}+1$. 

\item[-] Case $r_{i-1}=3$ and $r_{i}=3$. 

If $f(p^{+}_{i-1})=f(p_{i-1})=1$, Remark \ref{rem:inicios} implies that $f(p^{-}_{i})=f(p_{i})=1$. Then,  $f(P_{i-1})=2t_{i-1}+1$ and $f(P_{i})=2t_{i}+2$. 
Analogously, if $f(p^{-}_{i-1})=f(p_{i-1})=1$, Remark \ref{rem:inicios} implies that $f(p^{+}_{i})=f(p_{i})=1$. Then,  $f(P_{i-1})=2t_{i-1}+2$ and $f(P_{i})=2t_{i}+1$. 

In either case, for every pair of consecutive components of $R$ with value 3, we have 3 more units in the value of the optimal total dominating function $f$.

\end{itemize}
Let us now consider the value of this optimal total dominating function $f$ in \ref{eq: gamma prop} according to the value of $r_1$.

\begin{itemize}

\item If $r_1=2$, without loss of generality, we analyze the case $r_{i_0-1}=2$ and $r_{i_0}=3=\ldots=r_{{i_0}+\ell_j-1}=3$ for some $\ell_j$, for ${i_0}\geq 2$.

After Lemma \ref{segundoclaim} and 
Remark \ref{rem:inicios}, $f(p^-_{i_0})=f(p_{i_0})=1$. Therefore, the previous analysis shows that, for each pair of consecutive $r_{i-1}=r_i=3$ we have $f(P_{i-1})+f(P_i)=2 t_{i-1}+2t_i+3$. Obtaining $3\lfloor\frac{\ell_j}{2}\rfloor$ extra units in the value of  optimal total dominating function. 

When $\ell_j$ is odd, we also have that $f(p_{i_0+\ell_j-1}^-)=f(p_{i_0+\ell_j-1})=1$ and then $f(P_{i_0+\ell_j-1})=2t_{i_0+\ell_j-1}+2$.

\[\sum_{i=i_0}^{{i_0}+\ell_j-1} f(P_i)=\left\{
\begin{array}{ll}
   \sum\limits_{i=i_0}^{{i_0}+\ell_j-1} 2 t_i+ 3 \lfloor\frac{\ell_j}{2}\rfloor  & \text{ if } \ell_j \text{ even,} \\
   &\\
    \sum\limits_{i=i_0}^{{i_0}+\ell_j-1} 2 t_i+ 3 \lfloor\frac{\ell_j}{2}\rfloor +2 & \text{ if } \ell_j \text{ odd.} 
\end{array}
\right.
\]
which can be rewritten as, 
\[\sum_{i=i_0}^{{i_0}+\ell_j-1} f(P_i)=
\sum_{i=i_0}^{{i_0}+\ell_j-1} 2t_i+ 2\ell_j-\lfloor\frac{\ell_j}{2}\rfloor.\]

Therefore,
\[
\gamma_t(G)=f(P_1)+\sum_{i=2}^{s-1} f(P_i)+2s=2t_1+\sum_{i=2}^{s-1} (2t_i+E_i)+ \sum_{j=1}^m(2\ell_j-\lfloor\frac{\ell_j}{2}\rfloor) +2s.
\] 

\item If $r_1=3$ then $r_{\ell_1+1}=2$. Following a similar analysis as in the previous case, we have two different situations according to the parity of $\ell_1$. When $\ell_1$ is even, $f(p_{\ell_1+1}^+)=f(p_{\ell_1+1})=1$ and then $f(P_{\ell_1+1})=2t_{\ell_1+1}$. If not, $f(p_{\ell_1+1}^-)=f(p_{\ell_1+1})=1$ and then $f(P_{\ell_1+1})=2t_{\ell_1+1}+1$.  
Again, taking into account that for each pair of consecutive $r_{i-1}=r_i=3$ we have $f(P_{i-1})+f(P_i)=2 t_{i-1}+2t_i+3$,
\[\sum_{i=1}^{\ell_1} f(P_i)=\left\{
\begin{array}{ll}
   \sum\limits_{i=1}^{\ell_1} 2 t_i+ 3 \lfloor\frac{\ell_1}{2}\rfloor  & \text{ if } \ell_1 \text{ even,} \\
   &\\
    \sum\limits_{i=1}^{\ell_1} 2 t_i+ 3 \lfloor\frac{\ell_1}{2}\rfloor +1 & \text{ if } \ell_1 \text{ odd.} 
\end{array}
\right.
\]
which can be rewritten as, 
\[\sum_{i=1}^{\ell_1} f(P_i)=
\sum_{i=1}^{\ell_1} 2t_i+ \ell_1+\lfloor\frac{\ell_1}{2}\rfloor.\]
Also, the value of $f(P_{\ell_1+1})$ can be expressed as $f(P_{\ell_1+1})=2t_{\ell_1+1}+ \ell_1-2\lfloor\frac{\ell_1}{2}\rfloor$. Then
\[\sum_{i=1}^{\ell_1+1} f(P_i)=
\sum_{i=1}^{\ell_1+1} 2t_i+2\ell_1-\lfloor\frac{\ell_1}{2}\rfloor .\]

Proving that the result in \ref{eq: gamma reducido} follows.
\end{itemize}

Finally, if $I_2=\emptyset$ the previous analysis is simplified since we have $m=1$ and $\ell_1=s-1$. Then,   

\[
\gamma_t(G)=\sum_{i=1}^{s-1}2 t_i+2s+s-1+\lfloor\frac{s-1}{2}\rfloor =\sum_{i=1}^{s-1}2 t_i+3s-1+\lfloor\frac{s-1}{2}\rfloor.
\]
thus proving \ref{eq: gamma reducido sin I2}.
\end{proof}

Although the hypothesis in Proposition \ref{prop: G'} look quite restrictive, combining it with the reduction proposed in Definition \ref{def:reduced_graph}, this result is crucial for finding the value of the total domination number of a general $D$-caterpillar.

In the previous result, we did not consider $\D$-caterpillars with $I_0\neq \emptyset$. In this case, we split the graph $G$ into $\D$-caterpillar subgraphs with the property of having all the components of their $R$-vectors in their $(T,R)$-representation either equal to zero or all different from zero. More precisely, 

\begin{definition}\label{def: split}
Let $G$ be a $\D$-caterpillar with its $(T,R)$-representation such that $I_0\neq \emptyset$. We define $R^1=(r_1,\dots,r_{i_0})$ where $i_0-1\notin I_0$ and $i_0 \in I_0$. Then, for each $i$ such that $|\{i-1,i\}\cap I_0|=1$, we generate two vectors $R^{i-j}=(r_{i-j},\dots,r_{i-1})$ and $R^{i}=(r_{i},\dots,r_{i+k})$ with the additional property that 
either one of them is the zero vector, or it has no component at value zero. 
The vectors $T^i$ are defined accordingly, preserving their correspondence with $R^i$. We say that this is the \emph{$I_0$-split decomposition} of the $D$-caterpillar $G$ into the $\D$-caterpillars $G^j=G(T^j,R^j)$ for $j\in J$.
\end{definition}

Note that when $I_0=\emptyset$, the $I_0$-split decomposition of $G$ results in the original graph $G$.

\begin{example}\label{ex: split}
Consider the  $\D$-caterpillar $G$ with 12 support vertices defined in Example \ref{ex:12parents} with $(T,R)$-representation:
\begin{itemize}
    \item $T = (1, 1, 1, 0, 0, 1, 1, 0, 1, 1, 0)$
    \item $R = (0, 1, 1, 2, 3, 3, 0, 3, 3, 3, 1)$
\end{itemize}
In this case, $I_0=\{1,7\}$ and the $I_0$-split decomposition of $G$ is given by the  following vectors:
\begin{itemize}
\item $R^1 =(0)$, $R^2=(1, 1, 2, 3, 3)$, $R^7=(0)$ and $R^8=(3, 3, 3, 1)$.
\item $T^1 =(1)$, $T^2=(1, 1, 0, 0, 1)$, $T^7=(1)$ and $T^8=(0, 1, 1, 0)$.
\end{itemize}
Then the caterpillar $G$ gives raise to the $\D$-caterpillars subgraphs $G^j$ for $j\in \{1,2,7,8\}$.

\begin{figure}[h]
\centering
\includegraphics[scale=0.25]{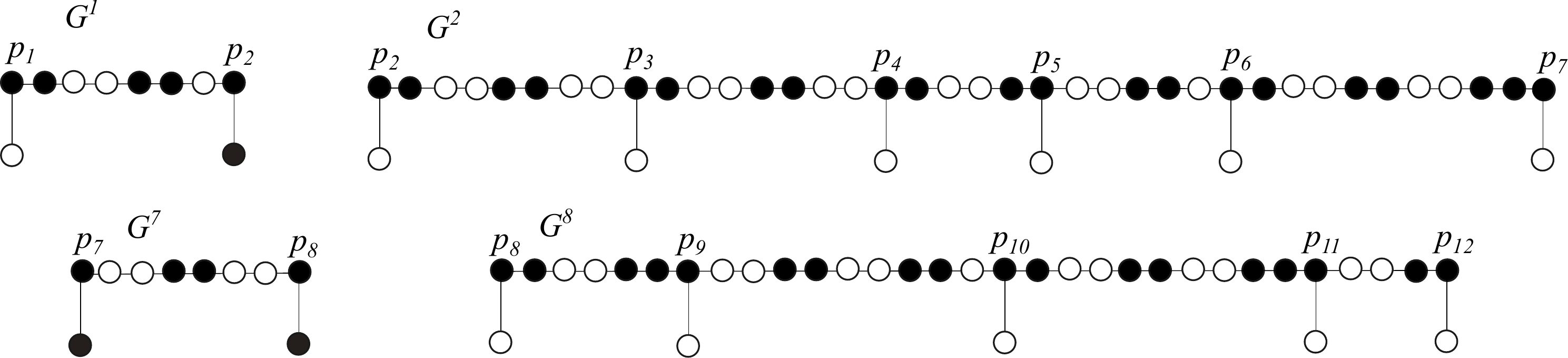}
\caption{The $I_0$-split decomposition of the $\D$-caterpillar $G$ in
Example \ref{ex:12parents}. The black vertices correspond to an optimal total dominating set in each  $G^{j}$ for $j\in \{1,2,7,8\}$.}
\label{fig: ejemI0_split}
\end{figure}

\end{example}


According to the definition of the $I_0$-split decomposition of a $\D$-caterpillar, we need to explicitly determine the value of the total domination number on two different type of $\D$-caterpillars. 

\begin{lemma}\label{lem: subcaterpillar}
Let $G$ be a $\D$-caterpillar with $(T,R)$-representation with $I_0=\emptyset$ or $I_0=[s-1]$.
\begin{itemize}
    \item If $I_0=\emptyset$ then 
    \begin{equation}\label{eq: I0 vacio}
     \gamma_t(G)=\gamma_t(G')+\sum_{i\in I_1}(2t_i+2)   
    \end{equation} 
    where $G'$ is the reduced graph obtained from $G$.
    \item If $I_0=[s-1]$ (i.e. $R$ is a zero-valued vector) then 
    \begin{equation}\label{eq: I0 novacio}
    \gamma_t(G)=\sum_{i=1}^{s-1} 2t_i+2s.
    \end{equation}
\end{itemize}
\end{lemma}
\begin{proof}
If $I_0=\emptyset$, after Proposition \ref{prop: G'}, we need to consider only the case when $r_i=1$, i.e., $i\in I_1$.
By Lemma \ref{primerclaim}, there exists an optimal total dominating function $f$ satisfying  $f(p_i^+)=f(p_i)=1$ or $f(p_i^-)=f(p_i)=1$. In any of these two possible situations, since $r_i=1$, Remark \ref{rem:inicios} shows that  $f(p_{i+1}^+)=f(p_{i+1})=1$ or $f(p_{i+1}^-)=f(p_{i+1})=1$, respectively.
Moreover, since $f(P_i)=2t_i$ and $f(N[p_i])=2$ we have that 
    \[\gamma_t(G)=\gamma_t(G')+\sum_{i\in I_1}(2t_i+2).\]

If $I_0=[s-1]$ the vector $R$ is the zero vector. We can assume that there exists an optimal total dominating function $f$ satisfying  $f(p_1^+)=f(p_1)=1$.  Since $r_1=0$, Remark \ref{rem:inicios} shows that  $f(p_{i}^+)=f(p_{i})=1$ for every $i\in [s-1]$ and then $f(P_i)=2t_i$. Finally, the value of $f$ in the closed neighbourhood of the support vertices proves the result. 
\end{proof}

We are now in a position to present a closed formula for total domination in a $D$-caterpillar, based on the result in Proposition \ref{prop: G'}, the reductions proposed in Definition \ref{def:reduced_graph}, Definition \ref{def: split} and Lemma \ref{lem: subcaterpillar}.

\begin{theorem}\label{thm: main}
Let $G$ be a $\D$-caterpillar with $(T,R)$-representation. The $I_0$-split decomposition of $G$  yields the collection $\{(T^j,R^j): j\in J\}$. For each $j\in J$, let $G^j=G(T^j,R^j)$. Then,
\[
\gamma_t(G)=\sum_{j\in J}\gamma_t(G^j) -2(|J|-1).
\]
\end{theorem}

\begin{proof}
The result clearly follows when $I_0=\emptyset$, due to Lemma \ref{lem: subcaterpillar} and Proposition \ref{prop: G'}.

Assume now that $I_0\neq \emptyset$ and let $\{(T^j,R^j): j\in J\}$ be its $I_0$-split decomposition.

Using the same reasoning as in the previous proofs, and taking into account the definition of each $G^j$, there exists an optimal total dominating function  $f$ in $G$ such that $f(p_j^+)=f(p_j)=1$ for each $p_j$ with $j\in J$. Then, in two consecutive graphs $G^{j_1}$ and $G^{j_2}$ (following the order of the vertices in the central path of $G$) we can sum-up the value of $f$ but taking into account its value on the closed neighbourhood of the support vertex they share in $G$. More precisely, if $G^{j_1}$ and $G^{j_2}$ are consecutive subgraphs in $G$ they share (in $G$) the last support vertex of $G^{j_1}$, $p_{j_1+1}$, and its adjacent leaf $h_{j_1+1}$.  Then, to obtain the value of the total domination number of $G$, we subtract two units (for the closed neighbourhood of the common support vertex) for each of pair of consecutive graphs in the $I_0$-split decomposition of $G$.
\end{proof}

\begin{example}
Let us now show how to apply the result in Theorem \ref{thm: main} for calculating the total domination number of the $\D$-caterpillar of Example \ref{ex:12parents}.

Since $I_0\neq \emptyset$, we perform the $I_0$-split decomposition, as shown in Example \ref{ex: split}.

In $G^2=G(T^2,R^2)$ the reduction in Definition \ref{def:reduced_graph} leads to $(G^2)'$ with $(R^2)'=(2,3,3)$ and $(T^2)'=(0,0,1)$. Similarly, $(G^8)'$ with $(R^8)'=(3,3,3)$ and $(T^8)'=(0,1,1)$.

In $G^1=G(T^1,R^1)$, by \ref{eq: I0 novacio}, we have 
\[\gamma_t(G^1)=2t_1+ 2s=2\cdot 1+4=6.\]

In $G^7=G(T^7,R^7)$, by \ref{eq: I0 novacio}, we have 
\[\gamma_t(G^7)=2t_7+ 2s=2\cdot 1+4=6.\]

 In Proposition \ref{prop: G'}, \ref{eq: gamma reducido} gives 
\[\gamma_t\left((G^2)'\right)=2 t_4+ 2t_5+2t_6 +2\ell_1-\lfloor\frac{\ell_1}{2}\rfloor+2.4=2.(0+0+1)+2.2-\lfloor\frac{2}{2}\rfloor+2.4=13\] and \ref{eq: gamma reducido sin I2} gives
\[\gamma_t\left((G^8)'\right)=2(t_8+t_9+t_{10})+\lfloor\frac{s-1}{2}\rfloor-1+3s=2(0+1+1)+1-1+12=16.\] 

Then \ref{eq: I0 vacio} implies that
\[\gamma_t(G^2)=\gamma_t\left((G^2)'\right)+2(t_2+t_3)+2s=13+2.(1+1)+2.2=21\] and \[\gamma_t(G^8)=\gamma_t\left((G^8)'\right)+2t_{11}+2=16+2.0+2=18.\]
Finally, since $J=\{1,2,7,8\}$,\[\gamma_t(G)=\gamma_t(G^1)+\gamma_t(G^2)+\gamma_t(G^7)+\gamma_t(G^8)-2(|J|-1)=6+21+6+18-6=45.\]

\end{example}
\section{Total $\{k\}$-domination in caterpillars} \label{sec: total dom general}

This last section is devoted to present the main result of this contribution, that is, a closed formula for the total $\{k\}$-domination in a general caterpillar graph. 

In Lemma \ref{interseccionnodisjunta} it is proved that we can recursively apply an extension of the original graph to obtain a caterpillar where all the support vertices have disjoint closed neighbourhoods, and, according to \ref{eq: gammat extG},  the $\{k\}$-domination number of the original graph can be obtained from the one of the extension.  Lemma \ref{lem: 1hijo} shows that without loss of generality we can assume that each support vertex has exactly one leaf adjacent to it.

Therefore, we can restrict to $\D$-caterpillar in the rest of the analysis. Finally, Theorem \ref{thm: 0-k} allow us to study the problem of total domination in $\D$-caterpillar, and gives the  $\{k\}$-domination number of the original graph in terms of the total domination number of it.
\begin{corollary}
Given a caterpillar $G$ with $r$ pairs of consecutive support vertices with non-disjoint closed neighbourhoods. Let $\Ext^r(G)$ the graph obtained from $G$ after performing the $r$ extensions defined in Lemma \ref{interseccionnodisjunta}. If $k$ is a positive integer number, then
 \[
 \gamma_t^{\{k\}}(G)=k \left(\gamma_t\left(\Ext^r(G) \right)-2r\right).
 \]   
\end{corollary}

\section{Conclusions and Future Work}
\label{sec: conclusions}

In this paper, we have investigated the total $\{k\}$-domination number, $\gamma_t^{\{k\}}(G)$, focusing on the class of caterpillar graphs. We established a general upper bound for this parameter and identified sufficient conditions for a graph to attain equality. Furthermore, we proved that this bound is tight for the family of caterpillars by showing that the total $\{k\}$-domination problem (T$\{k\}$DP) reduces to the classical total domination problem ($k=1$), yielding the relation $\gamma_t^{\{k\}}(G)=k\gamma_{t}(G)$. To achieve this, we introduced a novel $(T,R)$-representation for caterpillars based on the number of support vertices and the lengths of the paths induced between them along the central path. By exploiting this structural representation and analyzing the underlying congruence modulo 4, we obtained a closed-form expression for the exact total $\{k\}$-domination number of any caterpillar graph for all values of $k$.

Several interesting lines of research remain open. A natural extension of this work is to determine whether the reduction $\gamma_t^{\{k\}}(G)=k\gamma_t(G)$ holds for all trees, which would generalize our findings beyond the caterpillar family. Additionally, future research could focus on identifying other structural properties that enable wider classes of graphs to achieve the general upper bound provided in Proposition~\ref{prop1 y 2}, thereby expanding the contexts in which the T$\{k\}$DP can be simplified to the classical case.

\end{document}